\documentclass{amsart}

\usepackage{amsmath}
\usepackage{amssymb}
\usepackage{amsthm}
\usepackage{amsrefs}
\usepackage{amscd}
\usepackage{amsfonts}
\usepackage{graphicx}
\usepackage{pdfsync}
\usepackage{hyperref}
\usepackage{color}
\usepackage{verbatim}
\usepackage{enumerate}

\usepackage{caption}
\usepackage{subcaption}
\usepackage{array}
\usepackage{comment}

\newtheorem{theorem}{Theorem}

\theoremstyle{definition}

\newtheorem{remark}[theorem]{Remark}
\newtheorem{example}{Example}

\newtheorem{Def}[theorem]{Definition}

\newcommand{\inn}{\,\in\,}

\newcommand{\nn}{\newline\newline}

\def\R{{\mathbb R}}

\def\inn{\,\in\,}

\newcommand{\Bl}[1]{\textcolor{blue}{{#1}}}

\usepackage{tikz}
\usetikzlibrary{arrows,chains,matrix,positioning,scopes}
\makeatletter
\tikzset{join/.code=\tikzset{after node path={%
\ifx\tikzchainprevious\pgfutil@empty\else(\tikzchainprevious)%
edge[every join]#1(\tikzchaincurrent)\fi}}}
\makeatother
\tikzset{>=stealth',every on chain/.append style={join},
         every join/.style={->}}
\tikzstyle{labeled}=[execute at begin node=$\scriptstyle,
   execute at end node=$]

\begin{document}
\title{Modeling Collaborations with Persistent Homology}
\author[M.~Bampasidou]{Maria Bampasidou}
\address{Maria Bampasidou, FRE Department, University of Florida, 1179 MCCA, Gainesville, FL 32611, USA}
\email{mabampasidou@ufl.edu}

\author[T.~Gentimis]{Thanos Gentimis}
\address{Thanos Gentimis, ECE Department, NC State University, 2105D EBII, Raleigh, NC 27606, USA}
\email{agentim@ncsu.edu}

\begin{abstract}
In this paper we describe a model based on persistent homology that describes interactions between mathematicians in terms of collaborations. Some ideas from classical data analysis are used.
\end{abstract}

\date{}
\maketitle \tableofcontents

\section{Introduction}

The structure of social networks is one of the main objects of investigation in the study of social systems. Studies from the fields of psychology, mathematics \cites{Gro:95,Gro:02a,Gro:02b}, economics \cite{New:08}, biology \cite{Nic:11}, and computer science in an attempt to model associations among individuals create 1-dimensional graphs. The structure of those networks, one hopes, reveals the pattern of social interaction amongst agents (such as institutions, individuals, etc.).
An extensive literature exists for these graphs exploiting all possible angles such as combinatorics, graph theory, social sciences, statistics, random walks, adjacency matrix analysis, measure theory, complexity, logic, algorithmic computability and many more. Questions like: who is the most central member, how is the network affected by individuals, who is more influential and who is more peripheral, how can a network be improved and what are the reasons it fails, are thoroughly investigated in this framework.

In this paper, we focus on the network of scientific collaborations. Research relation systems and their study are not new to scientometrics. The field has a long history of the study of citation networks i.e. the network formed by the citations between papers. Also, co-authorship of a paper can be thought of as documenting a collaboration between two or more authors and these collaborations form what is now called a co-authorship network. Most of the times the networks mentioned above, assume that all agents are homogenous. We do not require that to be the case. We actually base our analysis on the differences amongst them.

We apply techniques from a new but very promising field of Computational Topology\cite{Car:09}. Persistent homology is used in an attempt to extract information from high-dimensional data sets recognising shapes the data provide us. We introduce the idea of a socioplex which captures ``connectivity" at higher levels (dimensions) and it is an improvement of the classic combinatorial sociograms since it contains them as it's one skeleton.
In terms of mathematical analysis, so far the goal is to create networks amongst mathematicians that are already collaborating. A lot of analysis has been done on co-authorship diagrams and relevant connections like the genealogy project for mathematicians, the co-authorship maps by the academic research society (created and maintained by Microsoft) and the Erdos number catalogue. Our model aims to predict future collaborations or at least indicate when two individuals should interact to produce something effectively.

We will not be using an already created network, for example co-authorship. Instead we are attempting to create a network of ``potential" collaborations amongst a discrete, finite number of agents, based on what we define as $R$-distance. This allows us to construct a multidimensional sequence of simplicial complexes indexed by a parameter $M$, corresponding to that distance. Part of our investigation is to find an effective $M$. In other words which is the most suitable $M>0$ such that if two agents are less than $M$ apart, according to our distance function, could collaborate.
Another aspect of our work is to identify persistent topological features, like connected components, holes, voids etc., and try to associate them with a physical meaning connected to a problem in the propagation of information (miss-collaborations).
We apply our model on the collaborations amongst mathematicians at Universities in the United States (set $X$) defining a distance function based on common characteristics of them.


The organization of this paper is as follows. In section 2, we recall the basic definitions from simplicial homology and persistent homology and establish the relative terminology. In section 3, we model our network by defining the metric to calculate distances amongst agents, and describe the program used to create our Socioplexes. Section 4, provides information about our data set and presents our findings in terms of persistent homology barcodes for the network of collaborations on mathematicians. We summarize our conclusions in section 5 and propose a generalization of the model to other fields.

We would like to thank Dr. Kevin Knudson, Dr. Ray Huffaker and Dr. Alexander Dranishnikov, for their very helpful conversations on these matters. All mistakes are our own.

\section{Preliminaries}

Consider a set of vertices $V$. A simplex of dimension $d$ is defined as the abstract combinatorial object $[v_1, v_2, ..., v_{d+1}]$. The geometric realization of a $d$-dimensional simplex, is the the convex hall generated by the $i-th$ standard basis vectors $\vec{e_i}$, $1\leq i\leq d+1$ in $\R^{d+1}$. A simplex of dimension 0, 1, 2, 3 is represented by a vertex, edge, triangle, pyramid etc., respectively.
Let $K$ be a finite, abstract, simplicial complex. We define $K$ as a union of simplices $\sigma_i$, $K=\sigma_1, \sigma_2, ..., \sigma_n$, such that for any two simplices $\sigma_i, \sigma_j$ we have either $\sigma_i\cap \sigma_j=\emptyset$ or $\sigma_i\cap\sigma_j=\sigma_k$, where $\sigma_k\inn K$.

The set of $i-$simplices, $\sigma^i$ (superscript denotes dimension), will be denoted  by $K^i$.
The number of simplices in a set $S$ is denoted by $|S|$. The dimension of the simplicial complex is $n\geq 0$ if $K^n\neq \emptyset$ and $K^{n+1}=\emptyset$.
Finite formal sums of simplices of $K^i$ with coefficients in a field (called $i-${\it chains}), define an additive abelian group structure on $K^i$.
In our case, the coefficients for those sums belong to $\mathbb R$ thus the group of $i-$chains, $C_i(K;\mathbb R)$ is a vector space over $\mathbb R$ with basis elements the simplices of dimension $i$.

If a simplex $\sigma$ is a face of another simplex $\sigma'$, we write  $\sigma\prec\sigma'$. We say that $\sigma'$ is a {\it coface} of $\sigma$.
A {\it proper face} of $\sigma\in K^i$, is a face of $\sigma$ of dimension $i-1$. The {\it boundary} of $\sigma$, denoted by $\partial(\sigma)$ is the formal sum (with coefficients in $\R$) of the proper faces of $\sigma$. The boundary operator is extended to all chains of $K$ by linearity.\par

An $i-$chain $a$ is an $i-${\it cycle} if $\partial_i(a)=0$ i.e. $a\inn\ker \partial_i$; it is an $i-${\it boundary}
if there is an $(i+1)-$chain $b$ such that $\partial_{i+1} (b)=a$ i.e. $a\inn \text{Im}(\partial_{i+1}) $. Two $i-$cycles $a$ and $a'$ are {\it homologous} if $a+a'$ is an $i-$boundary. It is not hard to show that $\ker \partial_i\subseteq \text{Im}(\partial_{i+1})$ \cite{Hat:02}
The quotient of $i-$cycles over $i-$boundaries is the  $i-$th  {\it homology space of} $K$  i.e.
$\displaystyle{H_i(K)=\frac{\text{Ker}(\partial_{i})}{\text{Im}(\partial_{i+1})}}$.

For any homology group $H_i(K)$ with finite dimension, we denote its dimension by $\beta_i=\text{dim}(H_i(K))$ and we refer to it as the {\it $i-$Betti number}. The basic topological structure of $K$ is quantified by the number of independent cycles in each homology space since the $i-th$ Betti number captures the number of independent $i-$dimensional surfaces. Particularly, $\beta_0$ represents the number of connected components of a space and $\beta_1$ counts the number of homological loops.

Given a discrete finite set of points $X$ in $\R^n$ and a finite distance $r$ one can built the Rips Complex, $K(X,r)$ which is the flag complex on the proximity graph of $X$. In other words if the nodes $v_1, v_2, ..., v_{d+1}$ form a set of diameter $r$ then they span the simplex $\sigma=[v_1, v_2, ..., v_{d+1}]$. For increasing values of $r$, one gets a nested sequence of simplicial complexes of $X$: \[ K_1\subseteq K_2 \subseteq ... \subseteq K_i \subseteq ... \subseteq K_k\]
Let $H_n^i=H_n(K_i)$, from the functorial properties of homology one gets a sequence of the form:
\[H_n^1\to H_n^2\to \dots H_n^i\to \dots\to H_n^k\]

\begin{figure}[ht!]
\centering
\includegraphics[scale=0.8]{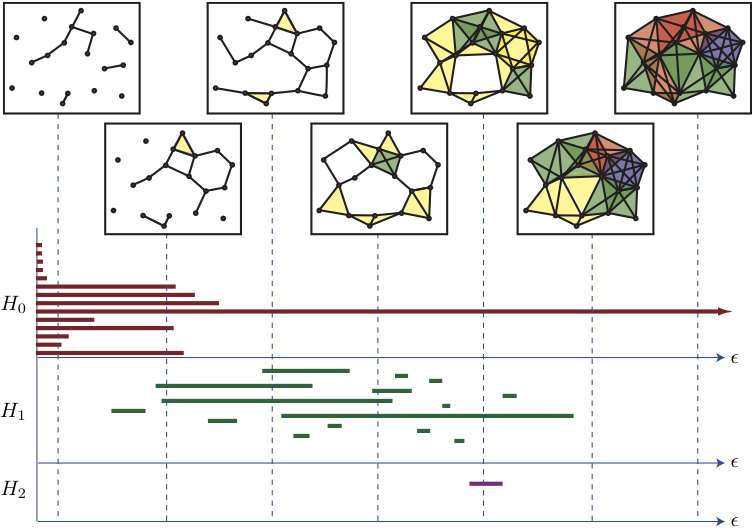}
\caption{An example of a filtration and the corresponding persistent barcodes}
\label{persistence}
\end{figure}

A particular class $[\alpha]$ may come into existence in $H_n(K_i)$ and then it either gets mapped to zero in $H_n(K_s)$ for some $s>i$ or to a non zero element in the last homology $H_n(K_k)$. This yields the persistent barcode, a collection of interval graphs lying above an axis parameterized by $r$. An interval of the form $[t, s]$ corresponds to a class that appears (is born) in $H^t_n$ and is mapped to zero (dies) at $H^s_n$. Figure \ref{persistence} gives a pictorial representation of the creation of these barcodes.

Classes that live to $H^k_n$ are represented by the infinite interval, to indicate that such classes are ``real" topological features of $X$. The collection of all barcodes (for all $n$) is the Persistent Homology diagram of the filtration corresponding to $X$ and its denoted by $Dgm(X)$.

\section{Distance Model}

From now on we define a set of agents $X$ which is a finite discrete set of elements. This can be made up of companies or individuals depending on the context. We see no reason why such a set could not contain both, but for simplicity we will later focus on the case where $X$ is made up of researches, specifically mathematicians.
We need to define a function that measures the research distance ($R-$distance) between two agents, A and B. Some of the key factors for this function could be: \begin{enumerate}[a)]
\item  the physical distance between two agents.
\item  the theoretical distance between the research projects of each agent.
\item  the willingness of each individual to trade research projects with another.
\item  the place where they received their schooling (especially if we are talking about researchers of a specific discipline).
\item  previous collaborations and previous citations.
\item  number of conferences they both attended, and more.
\end{enumerate}

Obviously the list is not exhaustive and one cannot know \textit{apriori} how much each of the aforementioned factors affects this theoretical distance. As a matter of fact some of them may even be untractable. We will thus focus on the factors that we can measure and try to estimate a good model by tweaking some weight parameters for these factors using a probabilistic method, Bayesian Estimation method.

The distance metric that we define is similar to the Hemming Distance for codes and it is a classic mathematical construction where one compares common (and uncommon) features between elements of a set. First we account for physical distance between two individuals. It is safe to say that people that live and work close by have higher chances to collaborate than people working in completely different parts of the world. Although that is not always true with today's modern communications (email, internet, teleconferences) physical distance seems to be a crucial factor for collaborations. So,

\begin{Def}Let $A$, and $B$ be two individuals. Then we define the physical distance between them $d_1(A, B)$ to be:
\begin{enumerate}[a)]
\item 0 if $A=B$,
\item 1 if $A, B $ work for the same institution,
\item 2 other
\end{enumerate}
\end{Def}

Second we account for their field of specialization. Assuming that we have a good handle on the pool of research interests these individuals have and we are comparing similar fields we can define the distance in research interests as follows:

\begin{Def} Let $A$, and $B$ be two individuals. Their distance in their research interests, $d_2(A, B)$ is defined as:
\begin{enumerate}[a)]
\item 0 if $A=B$,
\item 1 if $A, B$ have at least one similar subfield of research that they are working on,
\item 2 other
\end{enumerate}
\end{Def}
Notice that we use more than one fields of research for each individual in our comparison. For simplicity and easier computations in our example with mathematicians we used $3$ fields. The metric defined above is the simplest in this case. One could refine this metric by creating $3$ distances (one for each subfield) but we will use this simplified version in this paper. We will explore this idea in future papers.

From our experience a key factor of future collaborations between two individuals is the institution they graduated from. We consider people attending the same university, college more likely to collaborating. Lots of collaboration networks show the adviser as a central node and his students attached to him/her. Gradually, edges of collaborations sprout between the students themselves. This paper also is a collaboration between the authors, who attended and completed their PhD at the same University (UF). Therefore we define the distance in terms of ``schooling'' as follows:

\begin{Def} Let $A$, and $B$ be two individuals. Then:
\begin{enumerate}[a)]
\item $d_3(A, B)=0$ if they are the same person,
\item $d_3(A, B)=1$ if the two people graduated from the same university,
\item $d_3(A, B)=2$ other
\end{enumerate}
\end{Def}

It is clear that when two people already collaborate on a project it is very likely that they will collaborate again, especially in the case where the project is multifaceted or lengthy. Thus we consider the incidence of a past collaboration to be another factor that influences future collaborations so we define the relevant distance as follows:

\begin{Def} Let $A$, and $B$ be two individuals. Then:
\begin{enumerate}[a)]
\item $d_4(A, B)=0$ if they are the same person,
\item $d_4(A, B)=1$ if they have  collaborated in the past,
\item $d_4(A, B)=2$ if they haven't.

\end{enumerate}
\end{Def}

Finally from personal experience, especially in academia, researchers find potential collaborations between the people whose papers they have cited in their work, or who cite them in their papers. It is easier to work with like-minded people and people who have worked on similar questions as you have (and hence the citation connection). The relevant distance is defined as follows:

\begin{Def} Let $A$, and $B$ be two individuals. Then:
\begin{enumerate}[a)]
\item $d_5(A, B)=0$ if they are the same person,
\item $d_5(A, B)=1$ if one of them has cited the other in the past,
\item $d_5(A, B)=2$ if they haven't.

\end{enumerate}
\end{Def}

Of course one could compute various other distances, as we mentioned before. It is true that collaborations sometimes happen without any of the previous factors being in place. Still we believe that the following distance based on the previous 5 factors will wield a model that approximates nicely a ``closeness'' measure between mathematicians. We formalize our \textit{Research distance} in the following definition:

\begin{Def} Let $A$, and $B$ be two individuals. their $R$-distance is defined as: \begin{eqnarray*}d(A, B)&=&d_1(A, B)\cdot K_1+d_2(A, B)\cdot K_2+d_3(A, B)\cdot K_3\\&+&d_4(A, B)\cdot K_4+d_5(A, B)\cdot K_5\end{eqnarray*} where $K_1,K_2,K_3,K_4,K_5$ are the \textit{weights} for each factor.
\end{Def}

\begin{remark}  The weights $K_1,K_2,K_3,K_4,K_5$ are computed empirically using backtracking and Bayesian analysis. For that reason we require that \[\sum_{i=1}^5 K_i=1\]
\end{remark}

\begin{theorem} This $R$-distance is an actual metric on the set of all individuals we account in our sample $X$.
\end{theorem}
The proof is simple since all $d_i$'s are distances and the weighted sum of finitely many metrics is a metric.

\begin{Def}Let $M>0$. Define the $M$-neighborhood of $A$ to be the set of all people whose $R$-distance is smaller that $M$. We denote it by $N_M(A)$. We say then that $A$ and $B$ are $M$-close iff $B\inn N_M(A)$. This is nothing more than the definition of the $M$ ball around $A$ in the regular sense of metric spaces.
\end{Def}

\section{Mathematicians and M-Socioplexes}

In this section, we will use the mathematical formalism of persistent homology to infer topological information from a particular sample set on mathematicians in academics employed by U.S. public universities. We show how our distance model produces simplicial complexes associated to particular information about our sample. We then extract the Persistent Barcodes and give an interpretation based on observed data.

\subsection{Distances for Mathematicians}

We consider the following discrete distances: physical, mathematical interests, and specific indicators of past collaborations.
We obtained this information manually by looking at data-banks like $MathSciNet$ and the individuals web-pages.

With respect to the distance in mathematical interests we follow the classification of fields and subfields through $MathSciNet$, a database on publications through which mathematicians tend to characterize themselves. Obviously a mathematician can be part of many fields, but an easy way to classify them as such is to look at their official publications and collect the fields that they have published most of their work. For this analysis we looked at their 3 top fields. Information on collaborations, and citations is also documented in $MathSciNet$. Hence,

Given a set of mathematicians $X$  we define their $R$-distances as before: \begin{eqnarray*}d(A, B)&=& d_1(A, B)\cdot K_1+d_2(A, B)\cdot K_2+d_3(A, B)\cdot K_3\\ &+& d_4(A, B)\cdot K_4+d_5(A, B)\cdot K_5\end{eqnarray*}
The set $(X,d)$ is a metric space.

\subsection{Socioplexes}

Assume that one knows the distances between all elements in the set of mathematicians $X$. The following construction is a higher dimensional combinatorial analogue of the well known sociogram. Individuals are represented as nodes and their connections take the forms of edges, triangles, tetrahedra etc. The formal definition is the following.

\begin{Def} Fix a value for $M>0$. We define the $K_M$ socioplex as follows:
\begin{enumerate}[a)]
\item We assign a vertex for each mathematician.
\item We draw an edge between two mathematicians if one is in the $M$-neighborhood of the other.
\item We draw a 2-simplex (triangle) between three vertices A, B, and C if and only if each mathematician is in the $M$-neighborhood of the others.
\item We continue with higher dimensional simplices.
\end{enumerate}
\end{Def}

For different values of $M$ we obtain different socioplexes. It is obvious that for $M_1<M_2$ we have $K_{M_1}\subseteq K_{M_2}$.

\begin{example} An implementation of our algorithm for  threshold $60\%$ of the maximum distance between all the mathematicians at UF (diameter) yields the following M-Socioplex:
\end{example}
\begin{figure}
\begin{center}
\includegraphics[scale=0.4]{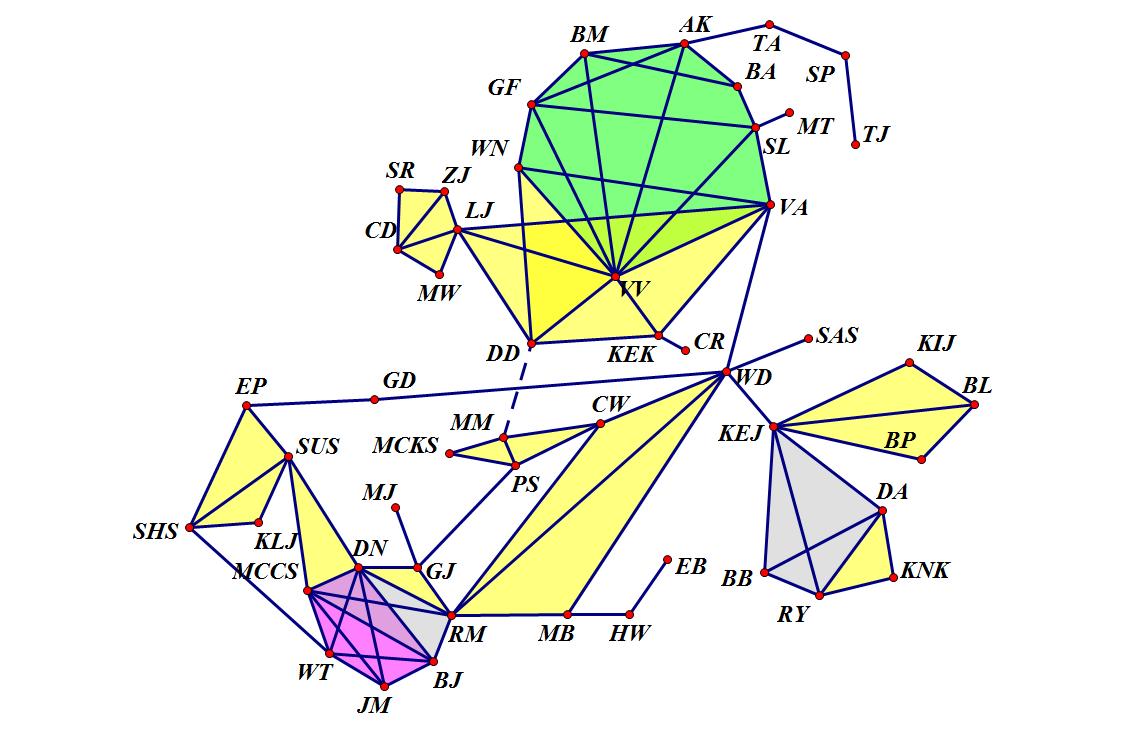}
\caption{The 6-Socioplex For the Mathematics Department At UF}
\end{center}
\end{figure}

An example of a creation of a Socioplex is given below:

\begin{example}[]

Suppose that we are given the following table of distances:
\[\begin{tabular}{r|c|c|c|c|c|c|c|c|c|c}
      & $u_1$  & $u_2$ &  $u_3$ &  $u_4$ &  $u_5$ & $u_6$ & $u_7$ &$ u_8$ & $u_9$ & $u_{10}$ \\ \hline
 $ u_1$ & 0    &  1  &  10  &  9   &  3   &  10 &  10 &  10 & 7   & 1\\
  $u_2$ &      &  0  &  1   &  2   &  10  &  10 &  10 &  10 & 10  & 10\\
  $u_3$ &      &     &  0   &  1   &  10  &  10 &  10 &  10 & 10  & 10 \\
  $u_4$ &      &     &      &  0   &  1   &  10 &  10 &  10 & 10  &10 \\
  $u_5$ &      &     &      &      &   0  &  1  &  4  &  3  & 5   &8 \\
 $ u_6$ &      &     &      &      &      &  0  &  1  &  2  & 10  & 10 \\
  $u_7 $&      &     &      &      &      &     &  0  &  1  & 10  & 10 \\
  $u_8$ &      &     &      &      &      &     &     &  0  & 1   &  10\\
  $u_9$ &      &     &      &      &      &     &     &     & 0   & 4\\
 $ u_{10}$&    &     &      &      &      &     &     &     &     & 0 \\
\end{tabular}\]
Then based on the threshold we choose we have the following Socioplexes \ref{socioplexes}:

\begin{figure}[ht!]
        \centering
        \begin{subfigure}[b]{0.2\textwidth}
                \includegraphics[width=\textwidth]{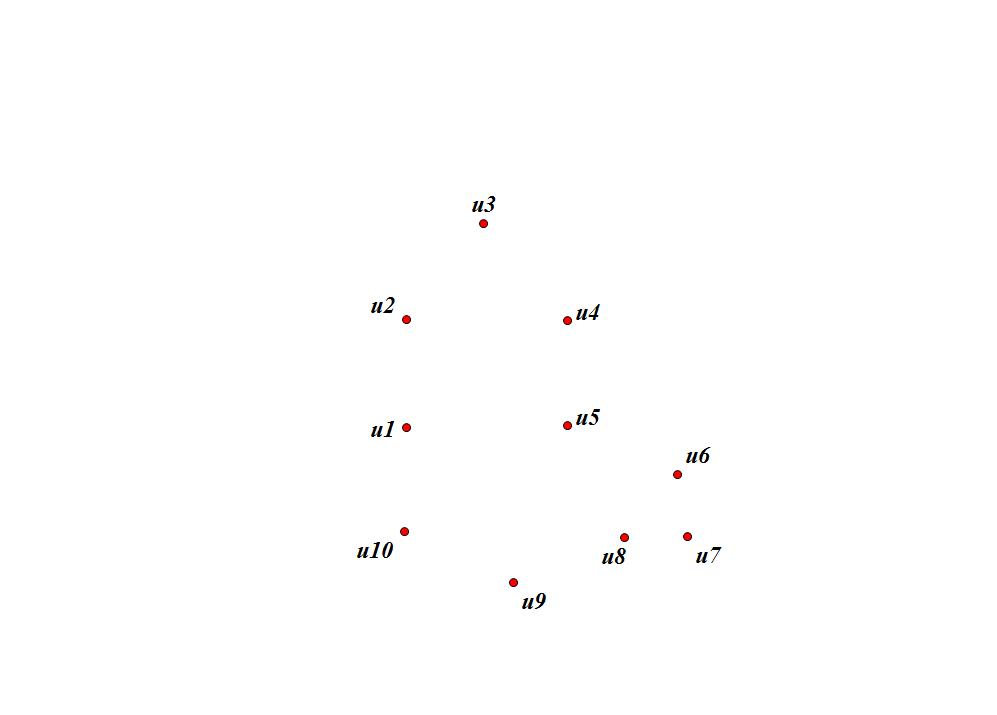}
                \caption{M=0}
        \end{subfigure}%
        \begin{subfigure}[b]{0.2\textwidth}
                \includegraphics[width=\textwidth]{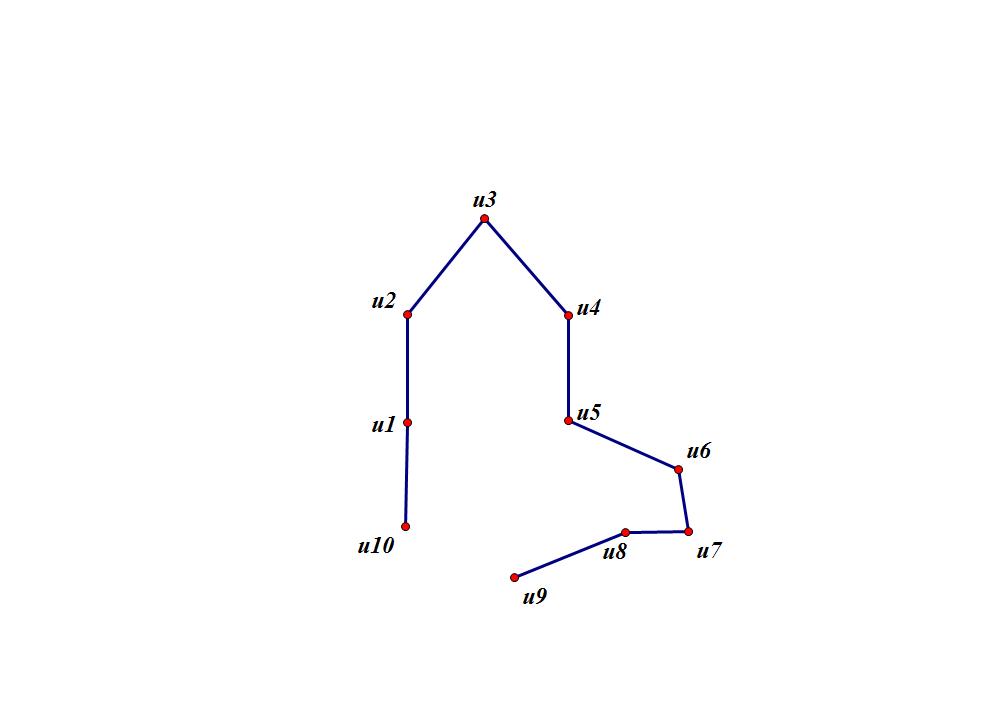}
                \caption{M=1}
        \end{subfigure}
        \begin{subfigure}[b]{0.2\textwidth}
                \includegraphics[width=\textwidth]{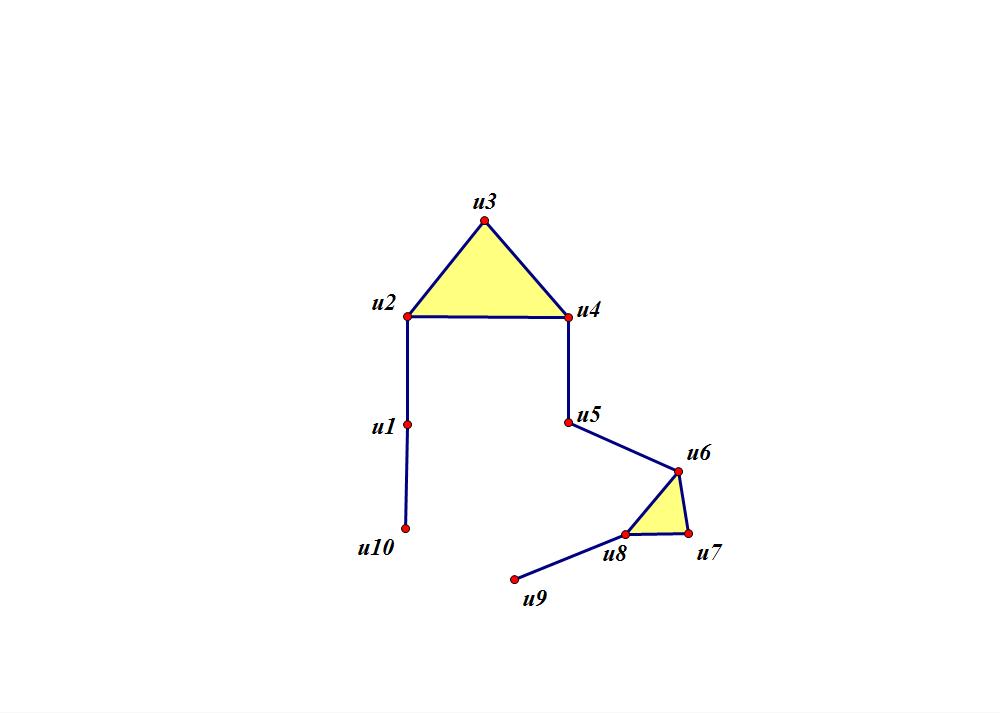}
                \caption{M=2}

        \end{subfigure}
         \begin{subfigure}[b]{0.2\textwidth}
                \includegraphics[width=\textwidth]{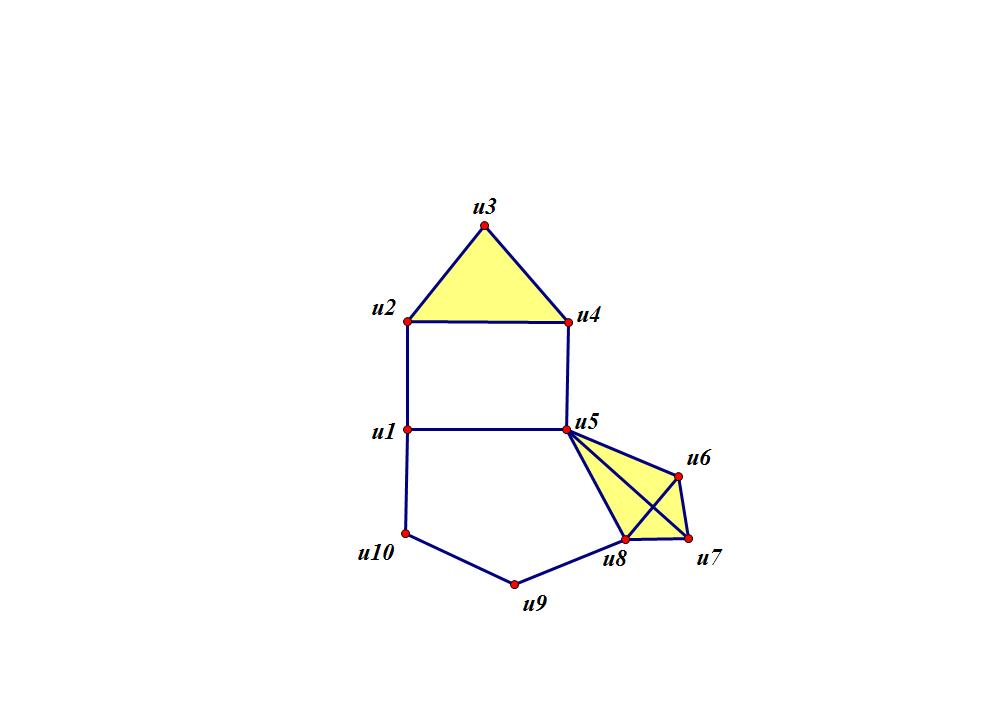}
                \caption{M=4}
                \label{M=4}
        \end{subfigure}

        \begin{subfigure}[b]{0.2\textwidth}
                \includegraphics[width=\textwidth]{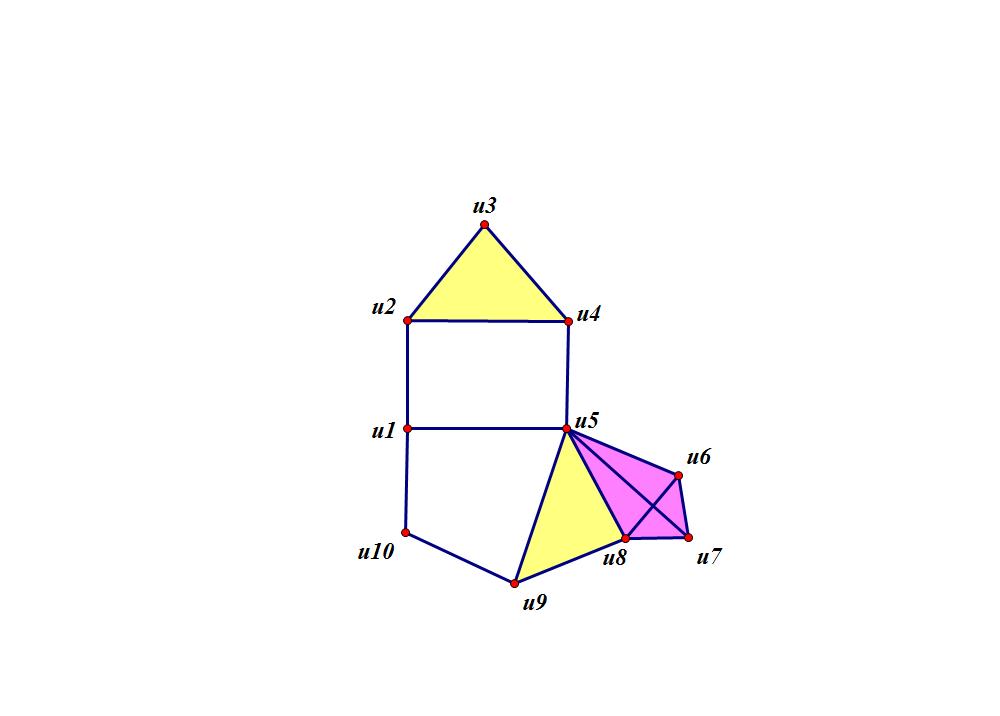}
                \caption{M=5}

        \end{subfigure}%
        ~ 
        \begin{subfigure}[b]{0.2\textwidth}
                \includegraphics[width=\textwidth]{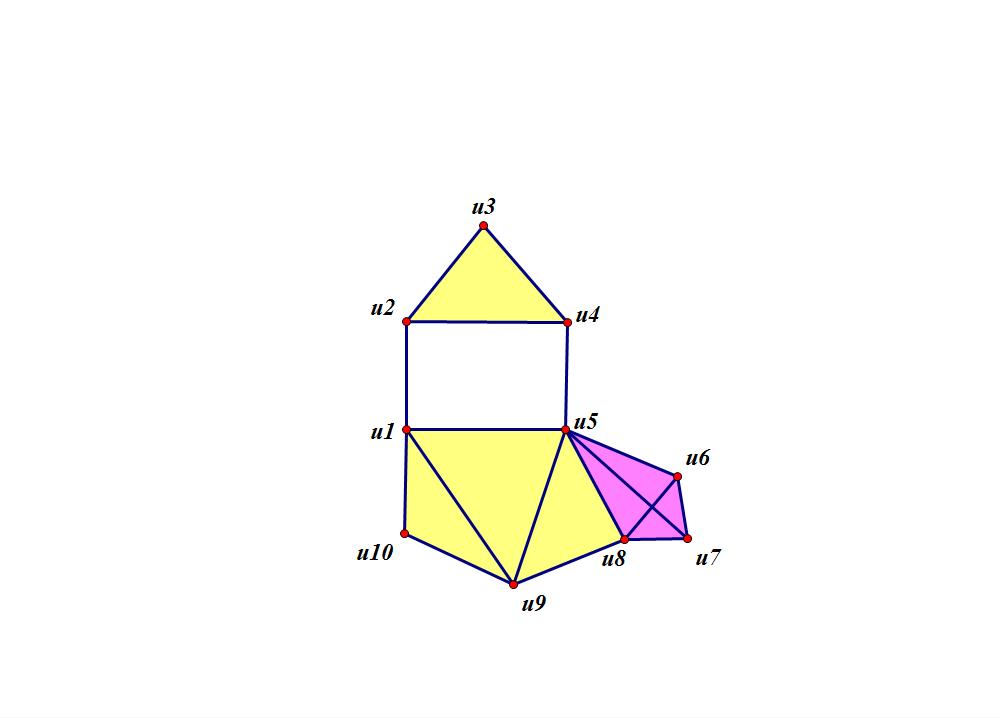}
                \caption{M=7}
        \end{subfigure}
        ~ 
        \begin{subfigure}[b]{0.2\textwidth}
                \includegraphics[width=\textwidth]{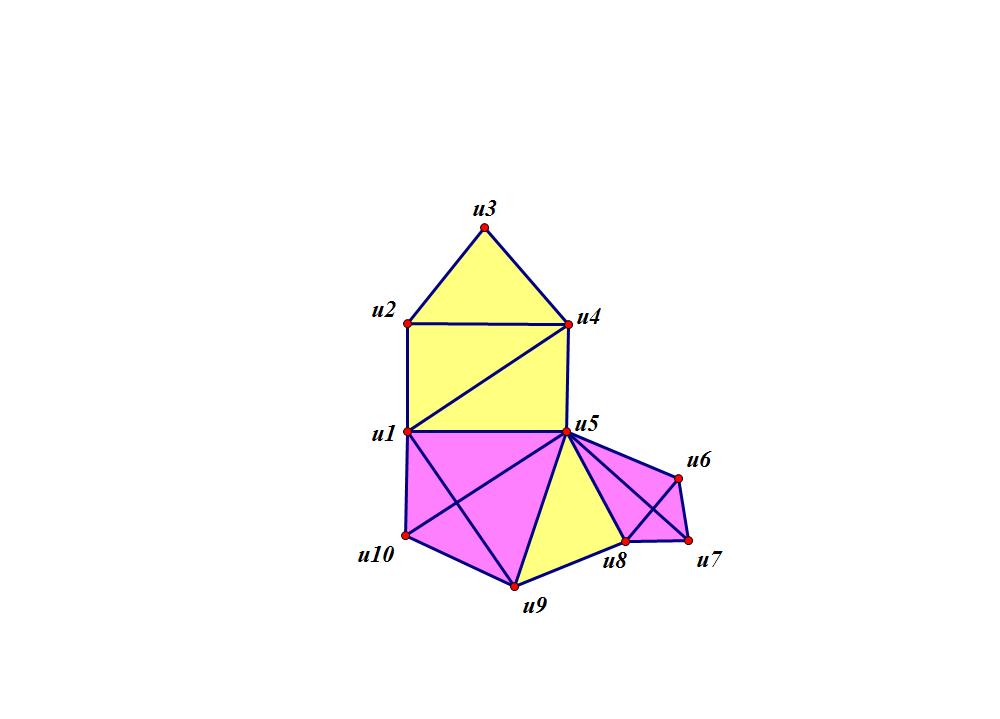}
                \caption{M=8}
        \end{subfigure}
        \caption{Socioplexes}
        \label{socioplexes}
\end{figure}
\end{example}

We claim that this higher dimensional pictorial representation contains more information than just the pair-wise interactions between mathematicians. The computations for this example where done by a simple matlab program described below.

The program takes as an input:
\begin{enumerate}[a)]
\item Professor's name.
\item Professor's location (University).
\item Professor's main fields (3 fields at most).
\item The institution where he/she attained his/her Ph.D.
\item Professor's list of citations (name-list).
\item Professor's list of collaborator (name-list).
\end{enumerate}

The program then generates the distances between any two individuals and given an $M$ creates the corresponding adjacency matrix which is the combinatorial structure needed to create the corresponding simplicial complex. Using available persistent homology software (javaplex and Perseus)\cites{Ada:11,Mis:13} we calculate the persistent homology and the corresponding barcodes for each M.

\begin{figure}
        \centering
        \begin{subfigure}[b]{0.9\textwidth}
                \includegraphics[width=\textwidth]{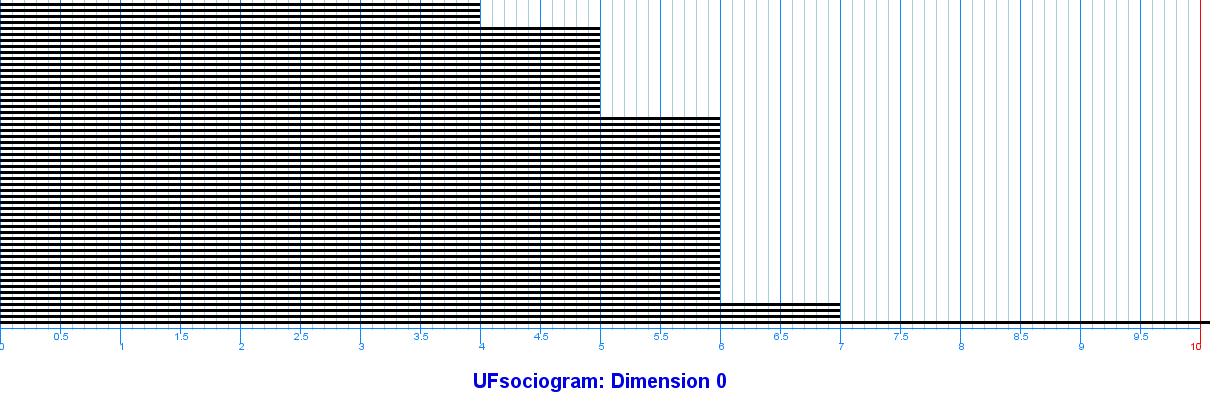}
                \caption{0-dimensional barcodes}
                \label{Uf0}
        \end{subfigure}%

        ~ 
        \begin{subfigure}[b]{0.9\textwidth}
                \includegraphics[width=\textwidth]{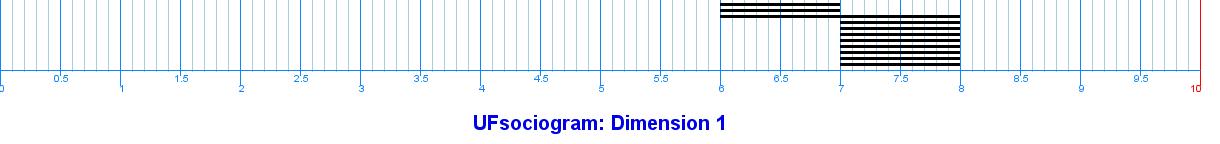}
                \caption{1-dimensional barcodes}
                \label{UF1}
        \end{subfigure}

        ~ 
        \begin{subfigure}[b]{0.9\textwidth}
                \includegraphics[width=\textwidth]{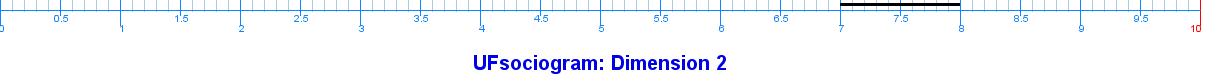}
                \caption{2-dimensional barcodes}
                \label{UF2}
        \end{subfigure}
        \caption{Persistent Barcodes for the mathematics department at UF}\label{barcodes}
\end{figure}

The corresponding persistent barcodes are depicted in figure \ref{barcodes}. By looking at those barcodes one can infer the following:
 \begin{enumerate}\item The $4$ longest zero bars at $70\%$ level correspond to the four big clusters:
  \begin{itemize}
  \item ``Applied mathematics'',
  \item ``Combinatorics and logic'',
  \item ``Analysis'', and
  \item  the rest.
   \end{itemize}
  \item The zero barcodes (clusters) at the $60\%$ level correspond very well to the natural division of the department into different fields.
  \item  There exist some short-lived $1$ dimensional bars in a region of the parameter which we call ``critical'' (more than one connected components and existence of loops, voids etc.).
  \end{enumerate}

\subsection {Obstructions to Collaborations}

The zero dimension persistent homology reveals the persistent clusters of research amongst mathematicians.

\begin{figure}[ht!]
\centering
\includegraphics[scale=0.8]{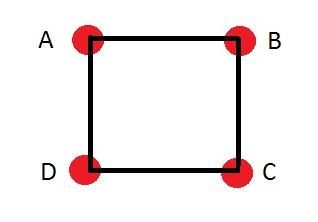}
\caption{1-dimensional holes as obstruction to Collaborations}
\label{hole}
\end{figure}

We propose that that the dimensional persistent homology (cycles), corresponds to ``obstructions to collaborations'' in the following sense. Consider the simplest case of four individuals $A,B,C,D$ which belong to the same cyclic ``research chain'' $A\leftrightarrow  B\leftrightarrow  C \leftrightarrow D\leftrightarrow A$ as depicted in figure \ref{hole}. The group that this four individuals create is not optimal. For instance individuals $B$ and $D$ don't communicate directly, which can be seen as a two hop information exchange. The ideal configuration would be the tetrahedron $[ABCD]$ where every two of them are connected or even any two triangles. Persistent long cycles would then reveal bigger ``weak'' information exchanges as the one we described above. Shorter cycles imply that the optimal configuration is easily obtainable by a small increase in the effective $M$. This corresponds to a small change in the research fields/interests of some of the individuals involved which theoretically may lead to a ``stronger group''.

The same idea could be extended to cycles of higher dimension, again somewhat quantifying the idea of ``weak'' and ``strong'' research teams.
One could argue that, given a good metric and effective parameter $M$ it suffices to calculate the standard homology of the $M$-Socioplex. But unlike persistent homology, classic homology is very unstable, i.e. it varies wildly with small perturbations of the metric (and therefore the distance between the points). Persistent homology on the other hand is a robust invariant of this model for the set of mathematicians, giving a concise description of the information flow dynamics and the potential of collaborations.

\section {Conclusions-Future work}

In this paper we presented a model based on characteristics of individuals to measure a theoretical distance between them in terms of future collaborations. Based on the model we created a graphic representation that lead to a network between individuals. We then used persistent homology to compute the corresponding barcodes.
We applied our analysis for the professors of the mathematics department at the University of Florida and gave an interpretation of the barcode output.

We claim that this model can be successfully used as a means to answer the following questions:
 \begin{enumerate}[a)]
 \item What can a university do to increase the interdepartmental collaborations?
 \item What can individual mathematicians do to increase their collaborations?
 \item What is a uniform measure of finding the strength of a grant proposal team?
\end{enumerate}
In the future we hope to be able to use the $MathSciNet$ data-bank to enlarge our set of mathematicians.
By getting access to a bigger data bank we will be able to answer the following questions:

\begin{enumerate}[a)]
 \item What are the ``correct'' coefficients $K_1, K_2, K_3, K_4, K_5$, meaning how much does the various distances (actual distance, distance in ideas and the willingness to share) affect the research distance.
\item What other features affect the ``distance'' between individuals, to what extend they affect (coefficients) and how can they be measured?
    \item How do mathematicians of different areas of mathematics collaborate (cluster)? Which Universities have more coherent research teams and is that reflected on their scoring?
 \end{enumerate}

\bibliographystyle{plain}
\bibliography{Bib_Ratopology}

\end{document}